\documentclass[11pt]{article}

\input{amssym.def}
\input{amssym}
\usepackage{amscd}

 \makeatletter\def\thechapter{\@Roman\c@chapter}\makeatother
\makeindex

\def\1{\hbox{\rm\rlap {1}\hskip .03in{\rm I}}}
\def\Int{{ \rm {Int}}}

\def\C{{  \Bbb C\,}}
 
\def\span{{  \rm {span}}}

 \def\Z {{\rm {Z}}} 
\def\skipaline{\vskip 12pt plus 1pt}

  \def\Hom{{ \rm  {Hom}}}
  
\def\diag{{ \rm  {diag}}}

  \def\exp {{ \rm  {exp}}}

  \def\dim  {{ \rm  {dim}}}

  \def\Tors {{ \rm  {Tors}}}
  \def\rk {{ \rm  {rk}}}

 \def\Ker {{ \rm  {Ker}}}
  \def\Im {{ \rm  {Im}}}
\def\diag {{ \rm  {diag}}}

\def\diag {{ \rm  {diag}}}  
\def\Z {{ \Bbb  {Z}}}
\def\Q{{ \Bbb  {Q}}}

\title{A homological estimate  for the Thurston norm}
\author{Vladimir Turaev}

\begin{document}
\maketitle

 {\bf Abstract} We establish a new 
 homological lower bound for the Thurston norm on
 1-cohomology of 3-manifolds. This  generalizes previous  results of
 C. McMullen, S. Harvey, and the author. We also establish an analogous lower
 bound   for 1-cohomology of 2-dimensional
   CW-complexes.

 {\bf AMS Classification} 57M27, 57M20, 57M05
 
 {\bf Keywords} 3-manifolds, 2-complexes, 
 skew fields,     Alexander 
  polynomial, Thurston  norm

\skipaline \skipaline 

 \centerline {\bf Introduction}
 
 \skipaline We establish a new 
 homological lower bound for the Thurston norm on
 1-cohomology of 3-manifolds. 
 This bound  generalizes a row of previously known results 
 and most notably the Seifert inequality saying that the genus of a knot in the 3-sphere
 is greater than or equal to half the span of the Alexander polynomial. 
 
 Recall first the definition of the Thurston norm, see [Th].  
 For a   finite   CW-space
$Y$, set 
$\chi_- 
(Y)=\sum_i \chi_- 
(Y_i) 
$ where 
the
sum runs over all connected components $Y_i$ of $Y$ and
$\chi_-(Y_i)=\max (-\chi 
(Y_i), 
0)$. 
 For a compact connected orientable 3-manifold $M$, the   
      Thurston norm  of $\psi\in H^1(M)=H^1(M;\Z)$ is defined 
by
 $\vert\vert  \psi
\vert\vert_T=\min_Y \, \chi_-(Y) $ where $Y$ runs
over  
 proper oriented
embedded (not necessarily connected) 
surfaces in $ M$   dual to $\psi$  with respect to some orientation of $M$.
 
  The recent interest to finding
  lower estimates for the Thurston norm on $H^1(M)$  was inspired by the
  famous (and much deeper) adjunction inequality in dimension 4, see [KM1] 
  and [MST]. 
  Suppose that $\partial M=\emptyset$. The adjunction 
  inequality applied to the 4-manifold $M\times S^1$ and to surfaces lying in the slice
  $M\times pt$ provides a finite set  
    $B\subset H_1(M)$ such
  that  $$\vert \vert \psi \vert \vert_T \geq \max_{h\in B} \vert \psi (h)\vert
  .\eqno (0.1)$$
  The set $B$ is defined in terms of spin$^c$-structures on
  $M$ and their Seiberg-Witten invariants. Namely, the 
  elements of $B$ are dual to the Chern classes of those 
   spin$^c$-structures  on
  $M$ whose Seiberg-Witten invariant is non-zero. 
  For more on this, see [Ak], [Au], [Kr],  
[KM2],   [Vi]. The Seiberg-Witten theory and the related
  Ozsvath-Szabo   theory yield  also 
  stronger estimates formulated in terms of Floer-type homology (see [OS]). 
   
   A homological approach  to  
  lower estimates for the Thurston norm was developed by C. McMullen [McM]
   who showed that the
  Thurston norm is bounded from below  by a norm derived from  the
  Alexander polynomial of $\pi_1(M)$. (Here we do not need to suppose that $\partial
  M=\emptyset$).  This result was generalized in two directions. 
  The author showed that
  the same is true for the twisted Alexander polynomials of $M$ associated with complex
  characters of the finite abelian group  $\Tors H_1(M)$, see [Tu1], [Tu2].
   For closed $M$, 
  the resulting set
  of inequalities is equivalent to  the adjuction inequality (0.1). 
  In an independent important development, 
  S. Harvey [Ha]
   derived   lower
  bounds for the Thurston norm  
    from   higher-order Alexander modules introduced
   in [COT] 
   (for  the exteriors of knots in $S^3$ the Harvey bounds were
    also obtained by T. Cochran [Co]).

  We now state our main result.
Fix a skew field $K$, i.e., a possibly non-commutative field with $1\neq 0$. 
Recall that any right (or left) $K$-module $H$ is free (see for instance [Coh1, Theorem
4.4.8]). 
The cardinality of a  basis of $H$ is 
 called the rank of $H$ and denoted $\rk_K H$.
It does not depend on the choice of the basis  (see [Coh1, Theorems
4.4.6-7]). If $H$ is finitely generated, then $\rk_K H$
is a non-negative integer.

Fix a field
automorphism $\alpha:K\to K$. One  defines the skew Laurent polynomial ring $\Lambda=
K[t^{\pm 1};\alpha]$
in one variable $t$ with coefficients in $K$ as follows.  The ring $ \Lambda$
consists of expressions $t^{-m} a_{-m} +...+t^{-1} a_{-1} +a_0+ ta_1+...+t^k a_k$ where
$m,k\geq 0$ and $a_i\in K$ for all $i$. These expressions add as ordinary Laurent
polynomials. Multiplication is defined on  monomials by $t^i a \,t^j b=t^{i+j} \alpha^j
(a) b$, where $i,j\in \Z$ and $a,b\in K$,  and then extended to arbitrary 
polynomials by additivity.  
 An easy computation shows that $\Lambda$ is an associative
ring with unit. Any $\Lambda$-module will be regarded as a $K$-module via the
inclusion $K\subset \Lambda$.

  For a pointed connected CW-space $X$, a 
  ring homomorphism 
 $ \varphi: \Z[\pi_1(X) ]$ $\to \Lambda $ allows us to regard $ \Lambda$    as a local system of
 coefficients on $X$ and to consider the corresponding right $ \Lambda$-module of 
 twisted homology $H_{\ast} 
 (X;\Lambda )$. (We shall discuss   twisted homology in   detail
 in Sect. 2.1). We say that $\varphi$ is {\it compatible with 
a cohomology class} $\psi\in  H^1(X)=H^1(X;\Z)$   if for any $\gamma\in \pi_1(X) $
 we have $\varphi (
 \gamma) \in  t^{\psi(\gamma)} K$ where $\psi(
 \gamma)\in \Z$ is the evaluation of $\psi
 $ on $\gamma$.

 \skipaline \noindent {\bf  Theorem 1.  } {\sl  Let
$  K, \alpha, \Lambda= K[t^{\pm 1};\alpha] $  be as above. 
 Let $M$ be a compact connected orientable 3-manifold. Set $\varepsilon_M=2$ if
$\partial M$ is   void or consists of 2-spheres  and set
$\varepsilon_M=1$ otherwise.
Let $\varphi : \Z[\pi_1(M)]\to \Lambda$ be a ring homomorphism 
  compatible with a non-zero cohomology class $\psi\in  H^1(M)$ and such that
   the right 
$\Lambda$-module 
 $H_1  (M;\Lambda ) $ is finitely generated over $K$. 
 If
the multiplicative group $\varphi(\pi_1(M))\subset \Lambda \backslash \{0\}$ is not cyclic, then
 $$\vert \vert \psi \vert \vert_T \geq \rk_K H_1  (M;\Lambda).\eqno
 (0.2)$$
 If
the group $\varphi(\pi_1(M))\subset \Lambda \backslash \{0\}$ is  cyclic,  then 
 $$\vert \vert \psi \vert \vert_T \geq \max (\rk_K  H_1  (M;\Lambda) -
  \varepsilon_M\,\vert \psi \vert ,0) \eqno
 (0.3)$$
where $\vert \psi \vert$ is the maximal positive integer dividing $\psi$ 
in the  
lattice $H^1(M)$. }

  \skipaline 
  Note that the number $\rk_K  H_1  (M;\Lambda)$ depends only on $\pi_1(M)$ and
  $\varphi$. 
  
   Theorem 1 generalizes   the  results of Seifert, McMullen, Harvey, Cochran, 
   and the author
   mentioned above. In
   particular,  Theorem 1 implies  the adjunction inequality in dimension 3. 
 
   The next theorem is a rather straightforward supplement to Theorem 1. 
   
    \skipaline \noindent {\bf  Theorem 2.  } {\sl  If under conditions of
    Theorem 1, $\psi$ is dual to the fiber of a fibration $M\to S^1$, then 
    $H_1  (M;\Lambda ) $ is finitely generated over $K$
    and (0.2), (0.3) become equalities. }

    \skipaline 
   The paper is organized as follows. In Sect. 1--3 we collected several 
   preliminary  definitions and lemmas used in Sect. 4 to prove Theorem 1. 
   Theorem 2 is proven in Sect. 5. 
   In Sect. 6 we
   discuss special cases and  relations with the previously known
   results. In Sect. 7 we consider  an   analogue of Theorem 1 for 2-dimensional
   CW-complexes. 
   
 \skipaline  {\bf  Acknowledgements.  } 
  I am   indebted to Peter Teichner for attracting my attention to the work of  
  T. Cochran and   S. Harvey. I am     indebted to T.
  Cochran and S. Harvey for sending me their preprints   [Co] and 
  [Ha].

  This paper was written during my stay at the Aarhus University. I am grateful
  to J. Dupont and Ib Madsen for hospitality. This stay was 
  supported by the European Union  TMR network \lq\lq Algebraic Lie
  Representations",  EC contract No ERB FMRX CT97-0100.
 
 \skipaline \skipaline \centerline {\bf 1. Codimension 1 subspaces of CW-spaces}

 \skipaline \noindent {\bf  1.1. Codimension 1 subspaces.} Let $X$ be a  
CW-space. A {\it codimension 1 subspace} of $X$ is a   CW-subspace $Y\subset X$ having  a
(closed) neighborhood   homeomorphic to $Y\times [-1,1 ]$ so that $Y=Y\times 0$.
It is understood that the product CW-structure on  $Y\times [-1, 1]$
should extend to a CW-structure of $X$ compatible with the given one (in the
sense  that they
 have a
common CW-subdivision). 

A typical example of a codimension 1 subspace is a proper two-sided $n$-dimensional PL-submanifold 
 of an $(n+1)$-dimensional PL-manifold. 
 (A submanifold is {\it proper} if it
 is compact and meets the boundary of the ambient manifold exactly along its own
 boundary).  The language of  codimension 1 subspaces is convient in our
 context. However only the case of proper surfaces in 3-manifolds will
 be actually needed for the proof of Theorem 1. The reader unwilling to
 consider the general setting of codimension 1 subspaces may safely restrict
 him/herself to this case.

 Each component $Y_i$ of a codimension 1 subspace $Y\subset X$ splits its cylinder neighborhood
into 2 components called the sides of $Y_i$. A choice of one of these sides 
(considered as   \lq\lq positive") determines a {\it coorientation} of $Y_i$.
We say that $Y$ is cooriented if all its components are cooriented.

  A cooriented codimension 1 subspace   $Y\subset X$ determines a 1-dimensional
  cohomology class  $\psi_Y\in H^1(X)$. Its value on a loop in $X$
  is   the algebraic intersection number of this loop with $Y$.  Here is
  a more formal definition of $\psi_Y$. Choose a cylinder neighborhood  
  $Y\times [-1, 1] \subset X$ of $Y=Y\times 0$ such that     $Y\times [0,1]$ lies
  on the positive side of $Y$. Define a map $g:X \to S^1=\{ z\in \Bbb C\, \vert\,
  \vert z\vert =1\}$
  which sends the complement of this neighborhood to $-1 $ and sends any point $
  (y,t)\in Y\times [-1,1]$ to $\exp (\pi i t)$. Then $\psi_Y=g^*(s)$ where $s$ is the
  generator of $H^1(S^1)=\Z$ determined by the counterclockwise orientation of $S^1$.
  
 We shall need a notion of a weighted codimension 1 subspace of $X$ which  
   formalizes subspaces with   
parallel
components.  A  codimension 1 subspace $Y$ of $X$ is {\it weighted} if each
  its component
$Y_i$ is endowed 
with an
 integer $w_i\geq 1$ called the {\sl weight}  of
$Y_i$. We  write
$Y=\cup_{i  } 
(Y_i, w_i)$ where $i$ runs over $\pi_0(Y)$.
 Such  $Y$  gives rise to an 
(unweighted)   codimension 1  
subspace $Y^{\#}\subset X$ obtained by replacing every 
  $Y_i$   
by
$w_i$ parallel copies $Y_i\times (1/w_i), Y_i\times (2/w_i),...,
 Y_i\times 1$ in a cylinder neighborhood of
$Y_i$.     A 
 coorientation of $Y$ induces a coorientation of $Y^{\#}$ in the
obvious way. 
For a cooriented  weighted codimension 1 subspace $Y=\cup_{i  } 
(Y_i, w_i)$ of $X$, set 
$$\psi_Y=\psi_{Y^{\#}}=\sum_{i }  w_i \, \psi_{Y_i} \in
H^1(X) 
.$$ 
If $Y$ is finite (as a CW-complex), then its components are  
 finite CW-complexes and there number is finite. Any $\Z$-valued 
 topological invariant
$\kappa$ of finite connected   CW-complexes extends to such $Y$ by linearity: 
$\hat \kappa 
(Y) =\sum_{i  }  w_i \,\kappa (Y_i) $.

\skipaline \noindent {\bf 1.2. Lemma.  } {\sl Let $\kappa$ be
a $\Z$-valued topological invariant
 of finite connected   CW-complexes. Let $X$ be a  
 connected CW-space and $Z\subset X$ be a   
 weighted  cooriented   codimension 1 finite CW-subspace.
  Then there is  a    weighted cooriented  
   codimension 1 finite CW-subspace $Y\subset X$
 such that $\psi_{Y}=\psi_Z$,  $\hat \kappa 
(Y)\leq \hat\kappa (Z)$,  and   $X\backslash
 Y$ is 
connected.}  

  \skipaline {\sl Proof.} 
We define several transformations on a  weighted cooriented  
codimension 1 subspace  
$Y=\cup_{i  
} 
(Y_i, w_i)$ 
of  $X$.  By \lq\lq decreasing the weight of  
$Y_i$ by 1" we mean
the transformation which reduces $w_i$ by $1$ and
keeps the other 
weights.  If
$w_i=1$, then this transformation removes $Y_i$
from $Y$.

Assume that $X\backslash   Y=X\backslash \cup_i Y_i  $ is not connected. For
a  component $N 
$ of 
$X\backslash    Y$, we define a {\it reduction} of 
$Y$ along   $N $. Let $\alpha_+$ (resp.\
$\alpha_-$)  be the set of 
all   
$i\in \pi_0(Y)$ such that $N$ is adjacent to $Y_i$ only on the positive (resp. negative)
side. 
Since $N\neq X\backslash    Y$ and $X$ is connected, at least one of the
sets 
$\alpha_+, \alpha_-$ is non-void.
Counting how many times   a loop in $X$ goes into or out of $N$, we
obtain  that 
$\sum_{i\in \alpha_{+}} \psi_{Y_i} =\sum_{i\in \alpha_{-}} \psi_{Y_i} $.
We modify $Y$ as follows.
   If $\alpha_{+}\neq \emptyset$ and $\sum_{i\in
\alpha_{+}}  \kappa 
(Y_i)\geq  \sum_{i\in \alpha_{-}}   \kappa
(Y_i)$, 
then we decrease by 1 the weights of all
$\{Y_i\}_{i\in 
\alpha_{+}}$ and 
increase by 1
the weights of   all $\{Y_i\}_{i\in \alpha_{-}}$.
 If $\alpha_{+}= 
\emptyset$ or  $\sum_{i\in 
\alpha_{+}}   \kappa (Y_i)< \sum_{i\in \alpha_{-}}
\kappa
(Y_i)$ 
 then we increase by 1 the weights of all
$\{Y_i\}_{i\in 
\alpha_{+}}$ and 
decrease by 1
the weights of   all $\{Y_i\}_{i\in \alpha_{-}}$.
   This
  yields another weighted cooriented  
   codimension 1 subspace  $Y'\subset X$ such that
$\psi_{Y'}=\psi_Y$ and $\hat \kappa (Y')\leq \hat \kappa
(Y)$.   
   Iterating 
this transformation, we eventually remove from
$Y$ at least one 
component
incident to $ N$ on one side.  We call this iteration
the   
 reduction 
of $Y$ along $N$.  The   reduction does not
increase $\hat \kappa $,
preserves $\psi_Y$ and 
strictly decreases the number of
components of $X\backslash   Y$.  
 If $  N$ is adjacent to only one component of
$Y$ 
 then the reduction 
removes 
  this component from $Y$.

We can now  prove the lemma. If $X\backslash   Z$
 is 
connected 
then  
$Z$ satisfies the requirements of the lemma.  If
$ X\backslash   Z$ is not 
connected 
then
iteratively applying to $Z$   
reductions 
along
the components of $ X\backslash   Z$ we eventually obtain
a weighted codimension 1 subspace $Y\subset X$
 such that $X\backslash   Y$ is connected.  Clearly, 
$\psi_Y=\psi_Z$ 
and  $\hat \kappa (Y)    \leq \hat \kappa (Z) $.

 \skipaline \skipaline \centerline {\bf 2. Twisted homology and the ring $\Lambda$}

  \skipaline \noindent {\bf  2.1. Twisted homology.} 
   Let $X$ be a connected CW-space
   with base point $x$. 
A  ring homomorphism 
 $ \varphi$ from $ \Z[\pi_1(X,x) ]  $ to an associative ring with unit $R$ 
 gives rise to a right $ R$-module of 
 twisted homology $H_{\ast} 
 (X;R )$  as follows. Consider the universal covering $p:\tilde X\to X$
 and denote by $\pi$ the group of covering transformations of $p$. 
To stay in line with [COT], [Co], [Ha], we consider $\tilde X$
as a {\it right} $\pi $-set. The right action of $\gamma\in \pi $
is just the usual left action of $\gamma^{-1}\in \pi $.
The cellular chain complex $C_\ast(\tilde X)$ of $\tilde X$ is then 
a right free chain complex over $ \Z[\pi ]$. To
 specify a basis of $C_\ast(\tilde
 X)$, it is enough to orient all  cells of $X$ 
 and choose their lifts   to $\tilde X$. 
 We can identify  $\pi=\pi_1(X,x)$ as usual and provide $R$ with the structure of   a
 $\Z[\pi ]$-bimodule via $arb= \varphi(a) r\varphi(b)$ for any $a,b  \in
 \Z[\pi ]$ and $r\in  R$. Then 
 $$H_{\ast}  (X;R )=H_{\ast}(  C_\ast (\tilde
 X)\otimes_{\Z[\pi]}  R ).  $$
 Note that the  $n$-th boundary
 homomorphism $C_n (\tilde X) \to C_{n-1} (\tilde X)$ is presented   by a matrix over 
 $\Z[\pi ]$. Applying $\varphi$ to its entries we obtain the matrix of the 
 $n$-th boundary
 homomorphism of the chain complex $C_\ast(\tilde
 X)\otimes_{\Z[\pi ]}  R$.

 The identification   $\pi=\pi_1(X,x)$ used above  depends on  a lift 
 of the base point $x\in X$ to
 $\tilde X$. 
For a different lift, this identification is composed with the
 conjugation by an element of $\pi$. 
 However,  the isomorphism class of the
 right $R$-module $H_{\ast}  (X;R)$ is easily seen to be well defined. 
  Similarly, the isomorphism class of   $H_{\ast}  (X;R)$ does not change when
  $\varphi$ is composed with conjugation by an invertible element of $R$.
 
   The
twisted homology extend  to   CW-subspaces $ 
A\subset X$ by 
$$H_{\ast}  (A;R)=H_*(  
C_{\ast}(p^{-1} (A)) \otimes_{\Z[\pi]} R). 
$$   Clearly, 
$H_{\ast}  (A;R)=\oplus_i H_{\ast}  (A_i;R)$ where $A_i$ runs over components of
$A$.

The twisted homology are invariant under cellular
subdivisions and form  
  usual
exact homology sequences such as the Mayer-Vietoris
  sequence.  Using a CW-decomposition
of $X$ with one 
$0$-cell,
one can check that $H_0(X;R)=R/  
\varphi (I) R  $ where $I$ is the augmentation ideal of $\Z[\pi_1(X,x)]$.   
 If $X$ is a closed   orientable $n$-dimensional PL-manifold, then using
 a CW-decomposition of $X$ with one $n$-cell  one easily computes that
 $H_n(X;R)=\{r\in R\vert\, \varphi(I)r=0\}$.

   A CW-subspace   $A\subset X$ is said to be 
   {\sl 
bad 
} (with respect to $\varphi$) if all loops in $A$   are mapped by 
$\varphi$ to $1\in R$.
This condition  does not depend on the way in which the origins of
the loops are   connected to the base point   $x\in X$.

\skipaline \noindent {\bf  2.2. The ring $\Lambda$.} We
  shall need a few   algebraic properties of the skew polynomial ring 
 $\Lambda=
K[t^{\pm 1};\alpha]$ defined in the introduction. 
Considering the highest degree terms of polynomials one immediately
obtains that $\Lambda$  has no zero-divisors. A slightly more
elaborated argument   shows that all right (and left) ideals
of $\Lambda$ are principal (see [Coh2, Prop. 2.1.1] or  [Co, Prop. 4.5]). This implies
that $\Lambda$ satisfies the right    Ore condition which says that  $a\Lambda\cap b
\Lambda\neq 0$ for all non-zero $a,b\in \Lambda$ (see [Coh2, Sect. 1.3]). Therefore the
formal expressions $ab^{-1}$ 
with $a,b \in \Lambda, b\neq 0$ form a skew field under  
natural addition and multiplication rules. 
 This field, called the classical right  field of quotients of
$\Lambda$, will be denoted  $Q$. The ring $\Lambda$ embeds into $Q$ via $a\mapsto
a1^{-1}$ so that $Q$ becomes a $\Lambda$-bimodule. 
It is known that $Q$ is a flat left $\Lambda$-module, i.e., the functor
$\otimes_\Lambda Q$ is exact (see [St, Prop. II.3.5]). 

For a right $\Lambda$-module
$H$,    the {\it torsion
submodule} $\Tors_\Lambda H$   of $H$ consists of 
  $h\in H$ such that $h a=0$ for a non-zero $ a\in \Lambda$.
Observe that if  $H$ is finitely generated over $K\subset \Lambda$, then
$\Tors_\Lambda H=H$. Indeed, for any $h\in H 
 $ the elements $h, ht, ht^2,...$ can not be  linearly  independent over $K$.
Therefore $x$ is annihilated by a non-zero element of   $  \Lambda$.
This   implies that for such $H$, we have $H\otimes_\Lambda Q=0$.

In Sect. 3 we shall use an estimate, established by Harvey [Ha],  on the $K$-rank of   $\Lambda$-modules
 presented by certain   matrices. Let $H$ be a right $\Lambda$-module with
 presentation matrix of the form $A+tB $ where $A$ and $B$ are  
  $(l\times m)$-matrices over  $K$ with $l,m=1,2,...$.
   Then 
  $\Tors_\Lambda H$ is finitely generated
   over
  $K$ and  $\rk_K \Tors_\Lambda H\leq \min (l,m)$, see [Ha], Prop. 9.1. (This is well
  known  when $K$ is a commutative field).
  
  \skipaline \noindent {\bf 2.3. Lemma.  } {\sl  Let 
 $X$ be a    connected CW-space with base point $x$. Let  
 $\varphi: \Z[\pi_1(X,x) ]\to \Lambda $  be a ring homomorphism 
 compatible with the   cohomology class 
 $\psi_Y\in  H^1(X)$ determined by a  weighted cooriented codimension 1 
  subspace 
   $Y\subset  X$ with finite number of components. Let 
  $c$ be the number of bad components of $Y$ and $d$
   be the number of bad components of $X\backslash Y$. Then
   $$c\leq d+ \rk_Q \,(H_1 (X; \Lambda)\otimes_\Lambda Q) \eqno (2.1)$$
   where $Q$ is the  right  field of quotients of
$\Lambda$.
 In particular, if     $H_1(X;\Lambda)$    has finite rank over
 $K$, then $c\leq d$.}

  \skipaline {\sl Proof.}  Note first that both $Y$ and $X\backslash Y$ have a
  finite number of components so that the numbers $c,d$ are well defined.

   We begin by computing the twisted 
  $0$-homology $H_0 (A;\Lambda)$ for a connected CW-subspace $A\subset X$
  such that  the restriction of  $\psi$ to $A$ is zero.    If $A$
  is bad then $H_* (A; \Lambda)$
is the usual untwisted homology of $A$ with  
coefficients in $\Lambda$ and in particular $H_0 (A; \Lambda)  
=\Lambda$. 
 Suppose that $A$ is not bad. Then $\varphi(\gamma)\neq 1$ for a certain $ \gamma\in 
 \pi_1(X,x) $
 whose conjugacy class is represented by a loop in $A$. 
 Since $\psi\vert_A=0$, we have
 $\varphi(\gamma)\in K$. Then $\varphi(\gamma)-1$ is 
   invertible in $K\subset \Lambda$.
 This implies that $H_0 (A; \Lambda) =0$. 
  
 Let $U=Y\times [-1,1]$ be a closed cylinder    
neighborhood of $  
Y=Y\times 0 $ in 
$X$ and 
$N=\overline {X\backslash U} $.
 The  Mayer-Vietoris   sequence of  the triple
$(X=N\cup U, N, 
U)$  gives 
an exact sequence
$$H_1 (X;\Lambda)\to H_0 (N\cap U;\Lambda)\to 
H_0 (N ;\Lambda) \oplus H_0 ( U;\Lambda).\eqno (2.2)$$ 
Observe that the restrictions of $\psi=\psi_Y$ to
 all components of $Y$ and  $N$ are
zero.  The computation 
 above shows that 
 $H_0 (N ; \Lambda)=  \Lambda^d$ and  
  $H_0 (U ; \Lambda) =\Lambda^{c}$. Similarly, $H_0 (N\cap U; \Lambda)=\Lambda^{2c}$. 
   Tensoring the exact sequence (2.2)
with   $Q$ we obtain an exact sequence 
$$H_1 (X; \Lambda)\otimes_\Lambda Q \to Q^{2c}\to Q^{c+d}.\eqno (2.3)$$ 
Note that     the $Q$-rank is additive with respect to short
exact sequences of $Q$-modules (all such sequences split). 
Therefore taking $Q$-ranks in (2.3) we obtain that 
$2c\leq c+d+ \rk_Q \,(H_1 (X; \Lambda)\otimes_\Lambda Q)$. This implies  (2.1).  

 \skipaline \skipaline \centerline {\bf 3. More on twisted homology}

\skipaline \noindent {\bf 3.1. 
Homology with coefficients in $K$.  } 
 Consider a   connected CW-space $X$ with base point $x$ and 
  a ring homomorphism 
 $\varphi: \Z[\pi_1(X,x) ]\to \Lambda $   compatible with a  cohomology class 
 $\psi \in  H^1(X)$. Let $A\subset X$ be a connected CW-subspace of $X$ such that
 the restriction $\psi\vert_A\in H^1(A)$ of $\psi$ to $A$ is zero.   
 Pick a base point $a\in A$ and   a path (or  a homotopy class of paths) 
 $\nu$ in $X$ leading from $x$ to $a$. 
 Pushing the origin of loops from $a$ to $x$ along $\nu^{-1}$ we obtain 
 an inclusion homomorphism $\pi_1(A,a) \to \pi_1(X,x)   $.  Extending it by linearity to
 the group rings and composing with $\varphi$ we  
 obtain a ring homomorphism $\Z[\pi_1(A,a)]\to  K\subset \Lambda$.
 This allows us to consider the corresponding twisted homology $H_\ast^\nu (A; K)$.
 To simplify notation, we shall sometimes 
   omit $\nu$.

 We define  a  homomorphism $g_A^\nu: H_\ast^\nu (A; K)\to 
  H_\ast (A; \Lambda)$  linear over the ring inclusion $K\subset \Lambda$. 
  Let $p_A:\tilde A\to A, p:\tilde X\to X$ be universal coverings with base
  points $\tilde a\in \tilde
  A, \tilde x\in \tilde
  X$ lying over  $a$ and $x$, respectively. Let $q=q_A:\tilde
  A\to \tilde X$ be the unique map such that $p q=p_A$ and 
   $\nu$ lifts to a
  path    in $\tilde X$ leading from $\tilde x$ to $ q(\tilde a)$. 
  The map $q$ and the inclusion $K\subset \Lambda$ induce  a chain homomorphism
  $$C_{\ast}(\tilde A) \otimes_{\Z[\pi_1(A,a)]} K \to 
  C_{\ast}(p^{-1} (A)) \otimes_{\Z[\pi_1(X,x)]} \Lambda. \eqno (3.1)$$
  The induced homomorphism in homology is $g_A^\nu:H_\ast^\nu (A; K)\to 
  H_\ast (A; \Lambda)$.
  
 Observe   that the chain $\Lambda$-complex on the right hand side of (3.1) is obtained from  the
 chain $K$-complex on the left hand side   via
  $\otimes_K \Lambda$  where $\Lambda$ is regarded as a left $K$-module via the inclusion 
  $K\subset \Lambda$. Clearly,   
$\Lambda$ is a free   $K$-module and therefore $\Lambda$ is flat over $K$. 
Therefore    $g_A^\nu$ induces an
  isomorphism
  $$ H_\ast^\nu (A; K) \otimes_K \Lambda=H_\ast (A; \Lambda). \eqno (3.2)$$
  The $K$-rank   of    $H_\ast^\nu (A; K)$ does not depend on
  the choice of   $\nu$. However, the
  image of $g_A^\nu$  does depend on $\nu$. An easy computation shows that for any 
   $\gamma\in \pi_1(X,x)$, we have 
    $\Im (g_A^{ \gamma\nu})=\Im (g_A^{\nu} ) \,\varphi(
   \gamma )$.
   
   The homology $H_\ast  (A; K)$ is natural with respect 
   to certain inclusions $A\subset
   A'$ where $A\subset A'$ are pointed connected CW-subspaces of $X$ such that
   $\psi\vert_{A'}=0$. Namely, let 
   $a\in A, a'\in A'$ be the base points of $A,A'$ 
   and let $\nu $ (resp. $\nu'$) be the
   distinguished path leading from $x$ to $a$ (resp. to $a'$). Suppose that the
   path $\nu^{-1} \nu'$ is homotopic (modulo
    the endpoints) to a path   in $A'$. Then    $q_A(\tilde A)\subset q_{A'} (\tilde {A'})
    $. This inclusion induces a  $K$-linear homomorphism 
    $H_\ast^\nu (A; K) \to 
 H_\ast^{\nu'} (A'; K)$ in the usual way.

 The homology $H_\ast  (A; K)$ can be extended to 
 non-connected CW-subspaces  $A\subset X$ with 
   $\psi\vert_A=0$ by  $H_\ast (A;K)=\oplus_i
  H_\ast (A_i;k)$ where $A_i$ runs over components of $A$. Here 
  we need to assume that all components of $A$ are pointed and $x$ is connected to
  their   base points
   by  distinguished paths. 
   
   If $A$ is   finite as a CW-space then $H_\ast  (A; K)$ is finitely generated
   over $K$. Moreover, $\rk_K H_n (A; K)$ does not exceed the minimal number of
   $n$-cells in a CW-decomposition of $A$. This follows from   definitions
    and the additivity of the rank with respect to short exact
   sequences.
   
 In particular, if $\psi=\psi_Y$ is 
  determined by a weighted cooriented codimension 1 subspace  
  $Y $ of $ X$, then $\psi\vert_Y=0$   and we can apply the constructions
  above to $A=Y $. By definition,  $H_\ast (Y;K)=\oplus_i
  H_\ast (Y_i;k)$. 
  
In the sequel, we call a  class   $\psi\in H^1(X)$
   {\it primitive} if   $\psi$  
  is not divisible by integers   $\geq 2$ 
in the lattice $H^1(X)$. 
   
 \skipaline \noindent {\bf 3.2. Lemma.  } {\sl  Let 
 $X$ be a finite   connected CW-space with base point  $x$ and let  
 $\varphi: \Z[\pi_1(X,x) ]\to \Lambda $  be a ring homomorphism 
 compatible with the  primitive cohomology class 
 $\psi_Y\in  H^1(X)$ determined by a  weighted cooriented codimension 1
   finite  CW-subspace 
   $Y=\cup_{i} (Y_i,w_i)$ of $ X$.
 Then for any $n\geq 0$, }
 $$\rk_K \Tors_\Lambda H_n (X; \Lambda) \leq \sum_i w_i\, \rk_K H_n(Y_i;K).\eqno (3.3)$$
 {\sl In particular, if   $H_n(X;\Lambda)$  has finite rank over
 $K$ then}
 $$\rk_K H_n (X; \Lambda) \leq \sum_i w_i\, \rk_K H_n(Y_i;K).\eqno (3.4)$$
 
 \skipaline {\sl Proof.} The right hand side of (3.3) does not
 change if  $Y$ is replaced by $Y^{\#}$. Therefore  we can assume that the weights of all
 components of $Y$ are equal to 1.  Then the right hand side of (3.3) is just
  $ \rk_K H_n (Y; K)$.
 
   Let $U=Y\times [-1,1]$ be a closed cylinder    
neighborhood of $  
Y=Y\times 0 $ in 
$X$ such that $Y\times [0,1]$ lies on the positive side of   $Y$. Set 
$Y_i^\pm =Y_i\times  (\pm 1)$ and  
$Y^\pm=\cup_i Y_i^\pm =Y\times  (\pm 1)$. Set $N=\overline {X\backslash U} $. 
 Clearly, $N\cap 
U= Y^+\cup  Y^-$. We  provide $X$ with  a CW-decomposition (compatible with
the given one) so that   $N, U, Y^+, Y^-$ are subcomplexes.

The  Mayer-Vietoris    sequence of  the triple
$(X=N\cup U, N, 
U)$ gives 
an exact sequence
$$\begin {CD}   
H_n (N;\Lambda)\oplus H_n (U;\Lambda) @>j>>
  H_n (X;\Lambda)  @>\partial>> H_{n-1} (N\cap U;\Lambda). \end {CD}  $$ 
 Set $J=\Im\, j=\Ker\, \partial \subset H_n (X;\Lambda)$. 
Note that $\Tors_\Lambda J=\Tors_\Lambda H_n (X;\Lambda)$.   Indeed 
  $H_\ast (N\cap U;\Lambda)= H_\ast (Y^-;\Lambda)
\oplus H_\ast (Y^+;\Lambda)$ and the right $\Lambda$-module $H_\ast (Y^{\pm} ;\Lambda)
  =H_\ast (Y^{\pm}; K) \otimes_K \Lambda$ is free since  $H_\ast (Y^{\pm}; K)$
   is a free
  $K$-module. Therefore $\Tors_\Lambda H_n (X;\Lambda)\subset \Ker\,
  \partial =J$ so that $\Tors_\Lambda H_n (X;\Lambda)=
  \Tors_\Lambda J$.   
  
    Writing down the next term
 of the Mayer-Vietoris sequence, we obtain   an exact sequence 
  $$  H_n (N\cap U;\Lambda)  \to 
H_n (N;\Lambda)\oplus H_n (U;\Lambda) \to J \to 0. \eqno (3.5)$$ Clearly, the inclusion homomorphisms 
  $ H_n (Y^{\pm};\Lambda)\to H_n (U;\Lambda)$ are   isomorphisms. 
  Therefore the  exact sequence (3.5) yields an exact sequence 
  $$\begin {CD} H_n (Y^+;\Lambda)   @>\eta >> H_n (N;\Lambda) @>>>
  J  @>>> 0 \end {CD} $$
 where   $\eta $
 can be computed as follows. We  identify $H_n (Y^+;\Lambda)=
 H_n (Y^{-};\Lambda) $ using the   
  homeomorphism $y\times 1\mapsto 
 y\times (-
 1):Y^+\to Y^{-}$ where $y$ runs over $Y$. Denote by   $\eta^{\pm} $
  the inclusion homomorphism 
 $H_n (Y^{\pm};\Lambda) \to H_n (N;\Lambda)$. 
 Then $\eta= \eta^+  - \eta^- $.

 Observe that the restrictions of $\psi$ to all components of $Y$ and $N$
 are zero. Since $\psi$ is primitive, there is   $\gamma\in \pi_1(X,x)$ such that
 $\psi(\gamma)=1$.   Let us 
 fix base points in all components of $N$ and choose   paths in $X$
 connecting them to   $x$. Composing if necessary these paths
 with loops   representing powers of $\gamma$, 
 we can assume that the algebraic number  of intersections of each of 
 these paths with $Y$ 
 is 0. Pick a  point 
 $y_i $ on
 every component $Y_i$ of $Y$ and fix the   base points $y_i\times 
 (\pm 1)$   
  on   
 $Y_i^\pm$. Let us connect the latter base points  
  to the base points of the adjacent components
 of $N$ (by paths in $N$) and then further to $x$ along the  paths chosen
 above. 
  Now we can consider the twisted homology $
 H_\ast (Y^\pm;K)=\oplus_i H_\ast (Y_i^\pm;K)$ and $H_\ast (N;K)
 =\oplus_j H_\ast (N_j;K)$ where $j$ numerates the components of $N$. 
 
  A basis of each $K$-module $
 H_n (Y_i^\pm;K)$   determines a basis of the $\Lambda$-module $
 H_n (Y_i^\pm;\Lambda)$   using   (3.2). Similarly, for each
 component $N_j$ of $N$ a basis of   $
 H_n (N_j;K)$   determines a basis of $
 H_n (N_j;\Lambda)$. With respect to these bases the   homomorphisms 
  $\eta^{\pm} $ are presented by matrices with entries in $K$. 
 The bases of $H_n (Y_i^+;\Lambda)=  H_n (Y_i^-;\Lambda)$
  obtained in this way from bases in $H_n (Y_i^+;K) , H_n (Y_i^-;K)$ 
  are   related by a matrix over $\Lambda$. We claim that   all its 
  entries lie in
   $tK$. Indeed, by the computations of Sect. 3.1,  
   these entries   lie in 
   $\varphi (\gamma_i)K$ where $\gamma_i$ is the loop formed by
 the distinguished path   leading from $x$ to $y_i\times (-1)$, the interval 
    $y_i\times [-1,1]$, and   the distinguished path leading from $y_i\times 1$ to
    $x$. The algebraic number of intersections of $\gamma_i$ with $Y$ is
 $1$. Therefore $\psi(\gamma_i)=1$ and  
 $\varphi (\gamma_i)\in tK$. 
  This proves that the homomorphism $\eta= \eta^+  -
  \eta^- 
  $ can be 
  presented by a matrix of the form $A+tB$ where $A,B$ are matrices over $K$. 
  The Harvey inequality stated 
  in Sect. 2.2   implies that  
  $\rk_K \Tors_\Lambda H_n (X;\Lambda)=
 \rk_K \Tors_\Lambda J   \leq \rk_K H_n (Y; K)$.
  
  Finally, if   $H_n(X;\Lambda)$  has finite rank over
 $K$, then   $\Tors_\Lambda H_n(X;\Lambda)$ $=$ $ H_n (X; \Lambda)$ so that (3.3) implies
 (3.4).

\skipaline \noindent {\bf 3.3. Remarks}. 
1.  The assumption that $\psi$ is primitive  
can be removed  using for instance the argument in [Th, p. 103], 
but we shall not need this. 

2. The right hand side of (3.3) can be
often easily estimated from above which may give useful estimates for  $\rk_K
\Tors_\Lambda H_n (X;\Lambda)$. 
Assume, for concreteness, that in Lemma 3.2, the space $X$ is a compact orientable 
3-manifold and $Y\subset X$ is a proper
cooriented  embedded surface    with weights of all components equal to
1. Then   $\rk_K H_1 (Y; K)\leq b_1(Y)$
 where $b_1$ is the first
Betti number. This follows from the fact that each
 component $Y_i$ either has  a
CW-decomposition  with $b_1(Y_i)$ one-cells (if $\partial Y_i=\emptyset$) or 
  collapses onto  a
graph with $b_1(Y_i)$ one-cells (if $\partial Y_i\neq\emptyset$). Another
useful estimate: 
 $$\rk_K H_1(Y;K)=-\chi(Y)+\rk_K H_0(Y;K)+\rk_K
 H_2(Y;K)$$
 $$=-\chi(Y)+c+c'\leq \vert \vert \psi_Y \vert\vert_T +c+c'$$
 where $c$ is the number of bad components of $Y$ and $c'$ is the number of
 closed bad components of $Y$.    For
 example, let $X$ be  the exterior of an oriented $m$-component link in $S^3$   
 and let  $\psi \in H^1(X)$ be the cohomology class taking values 
 $k_1,...,k_m\in \Bbb Z$ on the
 meridians of the link. One can easily show (cf. [Ha, Cor. 10.4])
 that   there is a Thurston norm 
 minimizing proper surface in $  X$ dual to $\psi$, having  
  $\sum_j \vert k_j\vert$ boundary circles and no closed
 components.
 Such a  surface has at most $\sum_j \vert k_j\vert$   components and therefore
  for  all $\varphi$ as in Lemma 3.2,
 $\rk_K \Tors_\Lambda H_1 (X; \Lambda)\leq \vert \vert \psi \vert\vert_T
 +\sum_j \vert k_j\vert $.

 \skipaline \skipaline \centerline {\bf 4. Proof of Theorem 1}

\skipaline \noindent {\bf 4.1.   } We first reduce the theorem to the case of 
primitive $\psi$. 
Any $\psi \in H^1(M)$ can be presented in the form   $\psi= n \psi'$  
 where  $n=\vert \psi \vert\geq 1$ and  $\psi' \in H^1(M)$ is primitive.  
  Consider the skew
 polynomial ring $\Lambda'=K[u^{\pm 1}; \alpha^n]$ in the variable $u$. There is a
  ring
 embedding $j:\Lambda'\to \Lambda$ which maps any monomial $u^m a$ to 
 $  t^{mn} a$  where
 $a\in K, m\in \Z$.  If 
a ring homomorphism 
 $\varphi: \Z[\pi_1(M)]\to  \Lambda$ is  compatible with $\psi$ then 
  it takes values  
in $j(\Lambda')$ and splits therefore as a composition of $j$ with a ring
homomorphism $\varphi':\Z[\pi_1(M)]\to \Lambda'$  compatible with 
$\psi'$. 
Observe that the right $\Lambda$-module $H_\ast  (M;\Lambda)$
viewed as a $\Lambda'$ module splits as a direct sum of $n$ copies of 
$H_\ast  (M;\Lambda')$ (this is true already on the level of chain complexes). 
Hence $\rk_K H_1 (M;\Lambda)= n \,\rk_K H_1  (M;\Lambda')$. This and the homogeneity of the
Thurston norm $\vert \vert \psi \vert \vert_T=
n \vert \vert \psi' \vert \vert_T$ imply that if 
 the claim of  Theorem 1
 holds for $\psi'$ 
then it holds also for $\psi$. 
We assume from now on that $\psi$ is primitive. 

\skipaline \noindent {\bf 4.2.  }  
 Fixing an orientation on $M$ we can identify orientations of surfaces in $M$ with their
 coorientations.  We can     represent $\psi\cap [M]\in H_2(M,\partial M)$ by a
proper oriented 
surface   $Z\subset M$  such that   $\chi_- 
(Z)=\vert\vert \psi \vert\vert_T$. In the notation of Sect. 1.1, we have 
 $ \psi_Z=\psi$.  We can regard  $Z$ as a
weighted codimension 1 subspace of $M$
with  weights of all components equal to $1$. 
We   apply   Lemma 1.2 to $\kappa=\chi_-$ and $X=M$.
  This lemma yields 
a weighted proper oriented 
surface   $Y=\cup_i (Y_i, w_i)$ embedded in $ M$
such that $M\backslash  Y$ is connected, 
$\psi_Y=\psi_Z=\psi$ 
and  $\hat \chi_- (Y)  \leq \hat \chi_- (Z)=\chi_- (Z)$.
The inequalities
$$\vert\vert \psi \vert\vert_T \leq \chi_- (Y^{\#})= \hat \chi_- (Y)\leq  \chi_- (Z)
=\vert\vert \psi \vert\vert_T$$
imply that  $ \vert\vert 
\psi  \vert\vert_T=\hat \chi_- (Y)=\sum_i w_i \,\chi_- (Y_i)$. Note that all
components $Y_i$ of $Y$ are proper embedded oriented surfaces in $M$.

We now separate 2 cases depending on whether $M\backslash Y$ is bad or not. 

 \skipaline \noindent {\bf 4.3.} Suppose
  that $
M \backslash Y$ is not bad with respect to $\varphi$. Since $H_1(M;\Lambda)$ is
finitely generated over $K$,  Lemma 2.3 implies that $Y$ has no bad components. 
For any component $Y_i$ of $Y$ we have   $ H_0 (Y_i; K)=0$, cf. Sect. 2.1 or the
proof of Lemma 2.3.  If $\partial Y_i\neq \emptyset$, then $Y_i$ is contractible onto a
1-dimensional complex and therefore $H_2  (Y_i; K)=0$.
If $\partial Y_i= \emptyset$ then using a CW-decomposition
of $Y_i$ with  one 2-cell,
one can easily check that $ H_2(Y_i; K)  =0$. Then $\chi(Y_i)=
\sum_{j\geq 0} (-1)^j \,\rk_K H_j (Y_i;K) =- \rk_K\, H_1 (Y_i; K) $ and 
$\chi_- (Y_i) =\rk_K\, H_1 (Y_i; K)$.  
By Lemma 3.2,
$$\rk_K H_1 (M; \Lambda) \leq \sum_i w_i\, \rk_K H_1(Y_i;K)=\sum_i w_i\, \chi_- (Y_i) 
=\vert\vert 
\psi  \vert\vert_T.$$
This implies the claim of the theorem in this case. 

\skipaline \noindent {\bf 4.4.} Suppose
  that  
$M \backslash Y$ is   bad with respect to $\varphi$. All components of $Y$ are parallel to surfaces lying in 
$M \backslash Y$ and therefore they all are bad. Since $H_1(M;\Lambda)$ is
finitely generated over $K$,  Lemma 2.3 implies that
$Y$ 
is 
connected. 
Since the dual cohomology class $\psi=\psi_Y$ is primitive, the
weight of  the 
only 
component of  
$Y$ is equal to $1$. 

It is clear that 
$\varphi(\pi_1(M))\subset \Lambda$ is a cyclic group generated by   
 the value of $\varphi$ on the homotopy class of a loop in $M$ 
 intersecting $Y$
once transversely. 
 We must therefore prove  the inequality
  (0.3). Since $\vert\vert\psi\vert\vert_T\geq 0$ and 
  $\vert \psi \vert=1$, we need to prove that 
  $\vert \vert \psi \vert \vert_T \geq  \rk_K  H_1  (M;\Lambda) -
  \varepsilon_M$.
  
   Set $b=b_1(Y)$. 
Assume first that $\partial M$ is   void or consists of 2-spheres so that 
$\varepsilon_M=2$.
Gluing 3-balls to the components of $\partial M$, we obtain a closed 3-manifolds with the
same 1-homology (twisted or untwisted) and the same Thurston norm. Therefore without loss
of generality we can assume that $\partial M=\emptyset$. Then $Y$ is a 
closed surface   and  $\vert\vert
\psi \vert\vert_T=\chi_- 
(Y)=\max 
(b- 2,0)$. Since $Y$ is bad the twisted homology
$H_1(Y;K)$ coincide  with the  untwisted homology, i.e.,  $H_1(Y;K)=K^{b}$.
By Lemma 3.2, $$\vert\vert
\psi \vert\vert_T \geq b-2 = 
\rk_K H_1(Y; K) -\varepsilon_M \geq \rk_K H_1 (M; \Lambda)- \varepsilon_M
.$$

Assume from now on  that $\partial M$ has a component of  non-zero genus so that 
$\varepsilon_M=1$. We shall prove below
that in this case $\partial Y\neq \emptyset$. Then $\vert\vert
\psi \vert\vert_T=\chi_- 
(Y)=\max 
(b -1,0)$.
As in the previous paragraph, $\rk_K H_1(Y;K)= {b}$. 
Hence, by Lemma 3.2, $b   \geq \rk_K H_1 (M; \Lambda)$.
Therefore $\vert\vert
\psi \vert\vert_T \geq b -1  \geq \rk_K H_1 (M; \Lambda) -\varepsilon_M$.

 It remains to prove that 
$\partial Y\neq \emptyset$. 
 We begin with algebraic preliminaries. 
   If the unit $1\in K$ is annihilated by a positive integer then  the minimal such
   integer is prime and 
    the additive subgroup of $K$ generated by $1$
is a  finite commutative  field. It is denoted $F$. If   $1\in K$ is not annihilated by 
 positive integers  then $\Q$ embeds in $K$ and we set $F=\Q\subset K$.  
  It is easy to check that   $F$ lies in the center of $K$ and any 
 automorphism of $   K$ acts on $F$ as
 the identity.  Let $\Lambda_0=F[t^{\pm 1}]$ be the usual
   commutative Laurent polynomial ring in variable $t$
  with coefficients in $F$. The  inclusion $F\subset K$ extends to a ring inclusion $
  \Lambda_0 \hookrightarrow \Lambda$.  It is clear that $K$ is a free  left   $F$-module and
  therefore $\Lambda$ is a free left $\Lambda_0$-module. Hence $\Lambda$ is
 a flat  left $\Lambda_0$-module.
 
  Consider  
the 
infinite
cyclic covering $\hat M\to 
M$    determined by $\psi$. Since $\psi$ is primitive, the manifold $\hat M$ is
connected. The 
surface 
$Y$ lifts to a 
homeomorphic surface $\hat Y
\subset \hat 
M$  splitting $\hat M$ into two connected pieces,
$  M_-$ and 
$ 
M_+$ such that $  M_-\cap  
M_+=\hat Y$.  Let $t$ be the generating deck
transformation of 
the 
covering 
$\hat  
M\to
M$ such that $t   M_+ \subset   M_+$. 
The action of $t$ on the  usual (untwisted) homology   $H_\ast
(\hat M;F)  $ makes it a  module over  
$\Lambda_0=F[t^{\pm 1}]$.
 It follows from   definitions and the flatteness of $\Lambda$ as a left $\Lambda_0$-module
  that   
 $$H_\ast (M;\Lambda)= H_\ast (\hat M; F) \otimes_{\Lambda_0} \Lambda =
 H_\ast (\hat M; F) \otimes_{F} K.$$  The assumption that $H_1 (M;\Lambda)$ has finite
 rank over $K$ implies   that $H_1 (\hat M; F)$ must be a finite dimensional vector
 space over $F$.
 
Suppose that  
  $ \partial
Y= 
\emptyset$ so that $Y\cap \partial M=\emptyset $.  Let $Z$ be a component of $\partial
M$ of genus $g\geq 1$.  Cutting $M$ open
along $Y$ we 
obtain 
a compact  connected  3-manifold homeomorphic to
$N=  
M_+\backslash 
t(\Int   
M_+)$. Therefore $\partial N$ consists of $\hat 
Y \cup 
t(\hat  Y)$, a copy of $Z$, and possibly other closed  surfaces.
For any $r\geq 1$, the  boundary of the 3-manifold
$N_r=N\cup tN\cup... 
\cup 
t^{r-1} N$ includes $\hat  Y \cup t^r(\hat 
Y)$ and $r$ 
copies of $Z$.
By a well known corollary of the Poincar\'e duality,
the image of the 
inclusion homomorphism 
$H_1(\partial N_r;F) \to H_1(  N_r;F)$ has
dimension $(1/2)\, 
\dim_{F} 
H_1(\partial N_r;F)\geq gr$. Therefore $\dim_{F}\,
H_1(N_r; F) 
 \geq gr$. On the other hand,   the kernel of the
inclusion 
homomorphism 
$H_1(N_r; F) \to H_1(\hat  M; F)$ has
dimension
$\leq 2\, \dim_{F}\,
H_1(\hat  Y; F) $.  This shows that   $
H_1(\hat  
M; F) $ can not be finite dimensional. This contradiction
implies that  
$ \partial Y\neq \emptyset$   and completes the
proof of the 
theorem. 
 
  \skipaline \skipaline \centerline {\bf 5. Proof of Theorem 2}

\skipaline  Suppose that $\psi$ is dual to the fiber, $Y$, 
of a fibration $M\to S^1$. Since $M$ is connected, all components of $Y$ are
homeomorphic to each other and represent the same class $\psi'=\psi/
\vert \psi\vert
\in H^1(M)$. The inequalities  (0.2) and
(0.3) for $\psi$ are obtained from the corresponding  inequalities for $\psi$
via multiplication by $ \vert \psi\vert$, cf. Sect. 4.1.
 Therefore, replacing  $\psi$ by $\psi' $ we can
assume that $\psi$ is primitive and $Y$ is connected.

If    $Y$ is  a 2-sphere  or  
 a 2-disc, then $\pi_1(M)=\Z$, $\vert \vert \psi \vert\vert_T=0$,
 and   (0.3) becomes an equality $0=0$.
   Assume from now on that    
 $Y$ is neither  a 2-sphere nor a 2-disc. Then  
    $\chi (Y) \leq 0$ and $ \chi_-(Y) =-\chi (Y)$.
     
 The universal covering $\tilde M\to M$ factors
through the infinite cyclic covering $Y\times \Bbb R\to M$ 
determined by $\psi$. Therefore $\tilde M=\hat Y \times \Bbb R$ where $\hat
Y\to Y$ is a regular covering of $Y$. This easily implies that $H_\ast
(M;\Lambda)= H_\ast (Y; K)$ as right $K$-modules. 
 In particular $H_1  (M;\Lambda ) $ is finitely generated over $K$.

If the fiber $Y\subset M$ is not bad with respect to $\varphi$, then 
$\varphi (\pi_1(M))$ is not a cyclic group and we must prove that (0.2) becomes an
equality. It suffices to prove the inequality opposite to (0.2). 
By the computations of Sect. 4.3,  
$$\vert \vert \psi \vert\vert_T  \leq  \chi_- (Y)=-\chi (Y)
=\rk_K H_1 (Y; K)=\rk_K H_1
(M;\Lambda).$$
Assume that   $Y\subset M$ is  bad with respect to $\varphi$. Then 
$\varphi (\pi_1(M))$ is   a cyclic group and we must prove that (0.3) becomes an
equality.  It suffices to prove the inequality opposite to (0.3). 
If $\partial Y\neq \emptyset$, then 
$$\vert \vert \psi \vert\vert_T \leq  \chi_- (Y) =-\chi (Y)=\rk_K H_1 (Y; K) -\rk_K H_0 (Y; K)$$
$$
=\rk_K H_1 (Y; K) - 1=\rk_K H_1
(M;\Lambda) - \varepsilon_M \vert \psi \vert$$
since $\vert \psi \vert=1$ and $ \varepsilon_M=1$. 
If $\partial Y= \emptyset$, then similarly
$$\vert \vert \psi \vert\vert_T  \leq -\chi (Y)=
\rk_K H_1 (Y; K) -\rk_K H_0 (Y; K)-\rk_K H_2 (Y; K)$$
$$
=\rk_K H_1 (Y; K) - 2=\rk_K H_1
(M;\Lambda) - \varepsilon_M \vert \psi \vert$$
since $\vert \psi \vert=1$ and $ \varepsilon_M=2$.

 \skipaline \skipaline \centerline {\bf 6. Special cases}
 
 \skipaline  Throughout this section $M$ is a compact connected orientable
 3-manifold with fundamental group $\pi$.

 \skipaline \noindent {\bf  6.1. Seifert's inequality.} Let $\Lambda=F[t^{\pm
 1}]$ be the usual commutative ring of Laurent polynomials 
 on the variable $t$ with  
 coefficients in a commutative field $F$. Pick $\psi \in H^1(M)$ and let 
 $\varphi$ be the ring homomorphism $\Z[\pi]\to \Lambda$
 defined by $\varphi(\gamma)=t^{\psi(\gamma)}$ for $\gamma\in \pi$. It follows from
 definitions that $\varphi$ is compatible with $\psi$ and 
 $H_1(M;\Lambda)=H_1(\hat M;F)$ where $\hat M\to M$ is the infinite cyclic
 covering determined by $\psi$. Clearly, $\varphi(\pi)\approx \Z$. 
 The inequality (0.3) says that if $H_1(\hat M;F)$ is finite dimensional over
 $F$ then 
 $$\vert \vert \psi\vert\vert_T \geq \dim_F \,H_1(\hat M;F) - \varepsilon_M \vert \psi\vert.
 \eqno (6.1)$$
 This was known to H. Seifert at least in the case where $F=\Q$, $M$ is
 the exterior of a knot $L$ in $S^3$ and $\psi$ is the generator of $H^1(M)=\Z$.
 In this case the condition $\dim_F \,H_1(\hat M;F)<\infty$ is always satisfied,
 the number $\dim_F\, H_1(\hat M;F)$ is equal to the span of the Alexander
 polynomial $\Delta_L$ of $L$, and $\vert \vert \psi\vert\vert_T=\max (2
 g(L)-1,0)$
 where $g(L)$ is the genus of $L$.  
 Clearly, $\varepsilon_M=\vert \psi\vert=1$. Therefore
 (6.1) yields $2 g(L)\geq \span \Delta_L$.

 One can modify a little the definition of $\varphi$.
  Namely, any group homomorphism $\sigma:H_1(M)\to
 F^*=F\backslash \{0\}$ gives rise to a ring homomorphism $ \Z[\pi]\to \Lambda$
 sending any $\gamma\in \pi$ to $ \sigma([\gamma]) t^{\psi(\gamma)}$ where   
  $[\gamma]\in H_1(M)$ is the homology class of $\gamma$. Applying Theorem 1 to
 this ring homomorphism we obtain   a \lq\lq twisted"
 version of  (6.1).

 \skipaline \noindent {\bf  6.2. McMullen's inequality.} McMullen [McM] defined
 for   any
 finitely generated group $\Gamma$ a   (semi-)norm
 $\vert\vert...\vert\vert_A$ on $H^1(\Gamma)$    called the Alexander norm.  
 It is derived from  the
 Alexander polynomial $\Delta_\Gamma$ as follows. Consider the free abelian
 group  $G=H_1(\Gamma)/\Tors
 H_1(\Gamma)$. Recall that $\Delta_\Gamma$ is an element of $\Z [G]$ 
 defined up to multiplication by $\pm 1$ and elements of $G$. Pick 
 a representative of $\Delta_\Gamma$ and  expand it as a finite sum 
 $ \sum_{g\in G} c_g g$ where $c_g\in \Z$. For any $\psi\in H^1(\Gamma)$,
 $$\vert\vert  \psi\vert\vert_A=\max_{g,g'\in G, c_g c_{g'}\neq 0} \vert \psi(g)
 -\psi(g')\vert $$
 where   $\psi (g)\in \Z$ is the
evaluation
of $\psi$ on $g$. 
 This norm does not dependent on the choice of the representative of $\Delta_\Gamma$.
  It is understood that if $\Delta_\Gamma=0$, then   
 $\vert\vert\psi\vert\vert_A=0$ for all $\psi$.
 
Applying these constructions to $\pi=\pi_1(M)$ we obtain  the Alexander norm
on $H^1(M)=H^1(\pi)$.
 McMullen [McM] proved that for all $\psi\in H^1(M) $,
  the Alexander   and  
 Thurston norms satisfy  $\vert \vert \psi \vert\vert_T \geq 
 \vert \vert \psi \vert\vert_A$ if $b_1(M)\geq 2$ and 
 $\vert \vert \psi \vert\vert_T \geq 
 \vert \vert \psi \vert\vert_A- \varepsilon_M \vert \psi\vert$ if $b_1(M)=1$.
 It is clear from Harvey's argument   [Ha, Prop. 5.12] that McMullen's inequality is
 a special case of Theorem 1. We shall prove this in a more general context in
 the next subsection.
  
 \skipaline \noindent {\bf  6.3.  Alexander-Fox norms.} 
 In generalization of McMullen's Alexander norm, the author [Tu1] defined 
 for   any
 finitely generated group $\Gamma$ a set of  
   (semi-)norms on $H^1(\Gamma)$
 numerated by $\sigma\in  \Hom (\Tors H_1(\Gamma), \C^*)$. 
 They are 
  called  Alexander-Fox norms  or (generalized) Alexander norms. Recall  
 their definition. Set $H=H_1(\Gamma), G=H/\Tors H$, and pick 
  $\sigma\in  \Hom (\Tors H, \C^*)$.
Fix a
splitting $H=\Tors H \times G$ (we use multiplicative notation for the
group operation in $H$). Consider the first elementary ideal 
$E(\Gamma)=E_1(\Gamma)\subset \Z[H]$
of $\Gamma$ (see [Fox]).  Consider
the ring
homomorphism $\tilde \sigma:\Bbb  Z[H] \to \C [G]$ sending $hg$
with $h\in
\Tors H, g\in G$ to $\sigma(h) g$ where $\sigma(h)\in \C^*\subset
\C$.
The ring $\C  [G ]$ is a unique factorization domain and we   set
$\Delta^\sigma (\Gamma)= \gcd  \, \tilde \sigma (E(\Gamma))$.  This $\gcd$
is an
element of $ \C [G]$ defined up to multiplication by elements of
$G$ and
nonzero complex numbers.     Pick
a
representative
$\sum_{g\in G } c_g g \in \C[G]$ of  $\Delta^\sigma(\Gamma)$.
For $\psi \in H^1(\Gamma)$, set $$\vert \vert \psi \vert \vert^\sigma=   
\max_{g,g'\in G, c_g  c_{g'} \neq 0}
\vert
\psi (g) - \psi (g') \vert . $$
 This (semi-)norm does not dependent on the choice of the representative 
 of $\Delta^\sigma (\Gamma)$ and does not depend on the choice of the splitting
 $H=\Tors H \times G$.
  If $ \Delta^\sigma (\Gamma)=0$, then the norm 
 $\vert\vert...\vert\vert^\sigma$ is identically zero.
   For $\sigma=1$, the norm $\vert\vert...\vert\vert^\sigma$ coincides with   
 McMullen's Alexander norm.
 
  The author proved in [Tu2] that for all $\psi\in H^1(M)=H^1(\pi) $ and all
  $\sigma\in  \Hom (\Tors H_1(M), \C^*)$,
    $$ \vert \vert \psi \vert\vert_T \geq \left \{ \begin {array} {ll}
   \vert \vert \psi \vert\vert^\sigma,~ {\rm {if}} 
\,\,\, b_1(M)\geq 2, 
\\
  \vert \vert \psi \vert\vert^\sigma- 
   \delta^1_\sigma\, \varepsilon_M\vert \psi\vert 
,~ 
{\rm
{if}} \,\,\,  b_1(M)=1.  \end {array} 
\right. \eqno (6.2)$$
Here $\delta^1_\sigma=1$ if $\sigma=1$ and $\delta^1_\sigma=0$
 otherwise. It is also explained in [Tu2] that varying $\sigma$ in (6.2)
  we obtain a set of inequalities   equivalent to
 the adjunction inequality (0.1).
 
 We now deduce  (6.2) from Theorem 1. Since all the norms
 involved in (6.2) are homogeneous, it suffices to prove (6.2) for primitive
 $\psi$.  Set $H=H_1(M)$
   and 
 denote by $U$ the kernel of the epimorphism $H/\Tors H\to \Z$ induced
 by $\psi$. Clearly,  $U$ is a free abelian group of rank $b_1(M)-1$. 
 The group ring $\C [U]$ is a commutative domain so that we can consider its
 field of quotients denoted $K$. Let $\Lambda=K [t^{\pm
 1}]$ be the   commutative ring of Laurent polynomials on the
  variable $t$ with  
 coefficients in $K$.
 
 Pick  $\tau \in H$ such that $\psi (\tau)=1$ and fix a
splitting $H=\Tors H \times U\times (\tau)$ where $  (\tau)$
 is the infinite cyclic
subgroup of $H$ generated by $\tau$.   
  Consider the ring homomorphism $\Z[H]\to \Lambda$
 sending $ hu\tau^m\in H$ to $\sigma(h) u t^{m}$ for any $h\in \Tors H, u\in
 U, m\in \Z$. Composing   with the obvious  ring projection
 $\Z[\pi] \to \Z[H]$ we obtain a ring homomorphism, $\varphi: \Z[\pi]  \to
  \Lambda$, compatible with $\psi$.
   
 \skipaline {\bf Claim.} {\it If $H_1(M;\Lambda)$  is finitely generated over
 $K $, then $\rk_K H_1(M;\Lambda)=\vert \vert \psi
 \vert\vert^\sigma$. If $H_1(M;\Lambda)$  is not finitely generated over
 $K $, then $\vert \vert \psi
 \vert\vert^\sigma=0$.}
 
 \skipaline This Claim and Theorem 1 imply (6.2). Indeed, if 
 $H_1(M;\Lambda)$  is not finitely generated over
 $K$ then   the right hand side of (6.2) is $\leq 0$ and (6.2)
 is obvious. Assume that $H_1(M;\Lambda)$  is   finitely generated over
 $K$. If $  \sigma=1$ and
 $b_1(M)=1$, then
 $\varphi(\pi )$ is a cyclic group and   (0.3) yields
 $$ \vert \vert \psi \vert\vert_T \geq \rk_K H_1(M;\Lambda)-
 \varepsilon_M \vert \psi \vert
 =\vert \vert \psi
 \vert\vert^\sigma- \delta^1_\sigma\,\varepsilon_M \vert \psi\vert.$$
 If $  \sigma\neq 1$ or
 $b_1(M)\geq 2$, then
 $\varphi(\pi )$ is not a cyclic group and   (0.2) similarly implies (6.2).
 
 Let us prove the Claim. For a  non-zero Laurent polynomial 
  $\lambda\in \Lambda$, we define the span $\span (\lambda)$ as the maximal
  difference $m-n$ where both $t^m$ and $t^n$ appear in $\lambda$  with
  non-zero coefficients.  Clearly, $\C [H/\Tors H]=\C[U] [\tau^{\pm 1}]$ 
  is the ring
  of Laurent polynomials on   $\tau$ with coefficients in $\C[U]$.
  Therefore   the inclusion $\C [U]\subset
 K$ extends to 
 a ring embedding $\nu:\C [H/\Tors H]\to \Lambda$ sending $\tau $
 to $t$. It follows   from definitions that 
 $\vert \vert \psi
 \vert\vert^\sigma=\span (\nu (\Delta))$ for any  representative $\Delta$ of 
 $\Delta^\sigma (\pi)$. Let $E\subset \Lambda$ be the ideal generated by 
 the set $\nu (\tilde \sigma (E(\pi))\subset \Lambda$. We claim that 
   $\nu (\Delta)=\gcd E$.   Recall that elements of a commutative ring 
 are {\it mutually prime} if  all their common divisors are invertible. 
  To prove the equality
   $\nu (\Delta)=\gcd E$ it suffices to show that  
    for any mutually prime elements of $\C [H/\Tors H]$,  
   their images in $\Lambda$ also are mutually prime. This easily follows from the
   fact that $\C [H/\Tors H] $ is a  unique factorization domain and the following
   property:  
  if  a non-zero element of $ \C[U]$ is divisible by $a\in \C [H/\Tors H]$
  in $ \C [H/\Tors H]=\C[U] [\tau^{\pm 1}]$, then
  $a\in \tau^m \C[U]$ for some $m\in \Z$. This property follows from the additivity of
  the span with respect to multiplication.

  The ideal $E\subset \Lambda$ is nothing but the first elementary ideal of the
  $\Lambda$-module $H_1(M,x; \Lambda)$ where $x$ is the base point of $M$. 
  This means that for any   presentation matrix of 
  $H_1(M,x; \Lambda)$ with $s$ rows and  $ \geq s$ columns,      $E $
   is generated by the $(s-1)\times
    (s-1)$-minors of this matrix.

    Consider   the exact homological 
sequence
$$0\to H_1( M; \Lambda) \to H_1( M,x; \Lambda) \to
H_0 (x; \Lambda) \to 
H_0 ( M; \Lambda).  $$ Clearly, $ H_0 (x; \Lambda)=\Lambda$ and
$H_0(M; \Lambda)=\Lambda/\varphi (I )\Lambda$ where
$I$ 
is the augmentation ideal of $ 
 \Z [\pi]$.  If $\sigma=1$ then 
$H_0 (M; \Lambda)=K$ and the inclusion homomorphism
$H_0(x; \Lambda) 
\to 
H_0  ( M; \Lambda) $ is the augmentation (summation of coefficients)
 $\Lambda=K[t^{\pm 1}]\to K$.  If $\sigma\neq
1$ then $H_0 (M; \Lambda)=0$.
In both cases the kernel   of the
inclusion 
homomorphism 
$H_0(x; \Lambda) \to H_0  ( M; \Lambda) $ is a free $\Lambda$-module of rank
one.  Hence $H_1 (
M,x; \Lambda)=H_1( M; \Lambda)\oplus \Lambda$.
   Therefore  for any   presentation matrix of the $\Lambda$-module
  $H_1(M; \Lambda)$ with $s$ rows and  $ \geq s$ columns,  the ideal  
   $E $ is generated by the $s\times
   s$-minors of this matrix. It remains to observe that $\Lambda$ is a
   principal ideal domain and therefore  $H_1(M; \Lambda)$ can be presented by
   a diagonal matrix $\diag (p_1,...,p_s)$ where $p_1,...,p_s\in \Lambda$.
  Then  $\nu (\Delta)= \gcd E=\prod_i p_i$ up to multiplication by non-zero elements
  of $K$.  If $p_i=0$ for some $i$ then $H_1(M;\Lambda)$  is not finitely generated over
 $K $ and $\vert \vert \psi
 \vert\vert^\sigma=  \span (\nu (\Delta))=0$. If $p_i\neq 0$ for all $i$ then $H_1(M;\Lambda)$  is   finitely generated over
 $K $ and has $K$-rank 
 $$\sum_i \span (p_i)=\span (\prod_i p_i)= \span (\nu (\Delta))=\vert \vert 
 \psi
 \vert\vert^\sigma .$$ 
 
  \skipaline \noindent {\bf  6.4. PTFA groups and Harvey's inequalities.} A
   source of skew fields  naturally
  appearing in the context of knot groups and 3-manifold groups, is the theory of 
   PTFA groups, cf. [COT], [Co], [Ha]. We call a group $\Gamma$    
  {\it poly-torsion-free-abelian} (PTFA) if
  it admits a normal series $\{1\}=\Gamma_0\subset \Gamma_1\subset...\subset \Gamma_n=\Gamma$ such
  that each   factor  $\Gamma_{i+1}/\Gamma_i$ is torsion-free abelian. Obviously,
  any PTFA
  group is torsion-free and solvable (the converse is not true).  
Any subgroup $H$ of a PTFA group $\Gamma$ is itself PTFA: for a 
  normal series $\{1\}=\Gamma_0\subset \Gamma_1\subset...\subset \Gamma_n=\Gamma$ as
  above, the groups   $H_i=\Gamma_i\cap H$ form a normal series
   $\{1\}=H_0\subset H_1\subset...\subset H_n=H$ such that  each factor $H_{i+1}/H_i$
   is a subgroup
   of $\Gamma_{i+1}/\Gamma_i$ and hence is torsion-free abelian. Any group $\Gamma$
   gives rise to a sequence of PTFA groups as follows. Consider the normal filtration
   $\Gamma=\Gamma^{(0)}\supset \Gamma^{(1)}\supset ...$ defined inductively by the
   condition that $ \Gamma^{(i)}/\Gamma^{(i+1)}$ is the maximal torsion-free abelian
   quotient of $\Gamma^{(i)}$ for $i=0,1,...$ It is clear   that $\Gamma/
   \Gamma^{(i)}$ is PTFA for all $i$. 
   
   For any PTFA group
   $\Gamma$ and any commutative field $F$, the group ring $F[\Gamma]$ is a right (and
   left) Ore domain,
   i.e., it has no zero divisors and satisfies the Ore condition mentioned in Sect. 2.2
   (see   [Pa,  pp. 588--592 and 611]). Therefore $F[\Gamma]$ embeds into its
   classical right ring of quotients, $F(G)$, which is a skew field. We shall use
   this construction to embed $F[\Gamma]$ into a skew Laurent
    polynomial 
    ring with coefficients in a skew field $K_\psi$ determined by 
 a primitive
  cohomology class
   $\psi\in H^1(\Gamma;\Z)$. 
   Namely, let
  $\Gamma_\psi$ be the kernel of the epimorphism $\Gamma\to \Z$ induced by $\psi$. Pick
  any   $t\in \Gamma$ such that $\psi(t)=1$. Every
  element of $\Gamma$ can be  uniquely expanded in the form $  t^m g$ with $m\in \Z$
  and $g\in \Gamma_\psi$. This allows us to identify 
   $F[\Gamma]$   with the skew Laurent polynomial ring 
  $F[\Gamma_\psi] [t^{\pm 1}]$   
    (the ring automorphism of $F[\Gamma_\psi]$ needed in the definition of the skew Laurent
    polynomial 
    ring is induced
  by the conjugation by $t\in \Gamma$). Therefore $F[\Gamma]\subset 
  K_\psi[t^{\pm 1}]$ where  
    $K_\psi=F(\Gamma_\psi)$. Note for
    the record that $ K_\psi [t^{\pm 1}]$ naturally embeds in $ F(\Gamma)$.

    Coming back to our geometric setting, we can use PTFA groups to construct 
   ring homomorphisms $\varphi$ as in Theorem 1.   Consider a PTFA group $\Gamma$
   and 
   a  commutative field $F$. Suppose that we have   group homomorphisms
   $\varphi_0:\pi=\pi_1(M)\to \Gamma$  and $\sigma:H_1(M)\to F^*= F\backslash
   \{0\}$.   
    Pick   a primitive cohomology class $ \psi\in H^1(\Gamma)$ and 
    recall the skew field $K_\psi$ and the embedding 
      $F [\Gamma]\subset K_\psi [t^{\pm 1}]$ constructed in the previous
    paragraph. Consider the ring homomorphism 
    $\varphi: \Z[\pi] \to K_\psi [t^{\pm 1}]$ which maps any $\gamma\in \pi$ to
     $\sigma ([\gamma]) \varphi_0(\gamma) \in 
     F [\Gamma]\subset K_\psi [t^{\pm 1}]$ where $[\gamma]\in H_1(M)$ is the
     homology class
     of $\gamma$. It is clear that 
    $\varphi$ is compatible with the
 cohomology class 
     $ \varphi_0^* (\psi)\in H^1(\pi)=H^1(M)$.  
  Thus we can apply Theorem
  1 to estimate $\vert \vert  \varphi_0^* (\psi)\vert\vert_T$ from
   $ \rk_{K_\psi} H_1(M; K_\psi [t^{\pm 1}])$. 
   
      Set now $F=\C$, $ \Gamma=\pi/
   \pi^{(i)}$ with $i\geq 1$, and let $\varphi_0$ be  the   projection  
  $\pi\to \Gamma$. For any $\sigma\in \Hom (H_1(M),\C^*)$ and 
 any primitive $\psi \in  H^1(M)$, we obtain an
 estimate  of $\vert\vert \psi\vert\vert_T$ (here we
     identify $H^1(\Gamma)=H^1(\pi)$ via $\varphi_0^* $). 
     This estimate extends by linearity to
     arbitrary   $\psi \in  H^1(M)$. In the
  case   $\sigma=1$, the
  resulting estimates were first obtained by S. Harvey [Ha]. 
  (She uses $F=\Q$ rather than $\C$ but this does not
    change the number $\rk_{K_\psi} H_1(M; K_\psi [t^{\pm 1}])$).  
    
     For $i=1$ we recover the estimates
 discussed in Sect.  6.3. 
  Fixing $i\geq 2$ and varying $\sigma$ in $\Hom (H_1(M),\C^*)$ we obtain  
  a  set
 of estimates which can be collectively considered as 
 the adjunction inequality of level $i$.
 This set of estimates is actually finite: one can check that  the number 
     $\rk_{K_\psi} H_1(M; K_\psi [t^{\pm 1}])$ does not change when 
     $\sigma$ is mulitplied by a homomorphism $H_1(M)\to \C^*$ sending $\Tors
     H_1(M)$ to $1$.

 Note finally that in this construction 
 $\varphi(\pi)\approx \Z $  if and only if $b_1(M)=1$,
 $\sigma (\Tors H_1(M))=1$,  
 and   either   $i=1$   or $i\geq 2$ and 
 $\pi^{(1)}=\pi^{(2)}=\pi^{(3)}=...$ In the case where $i\geq 2$ and 
 $\pi^{(1)}=\pi^{(2)}=\pi^{(3)}=...$ the inequality (0.3) is of no
 interest since  then $H_1(M; K_\psi [t^{\pm 1}])=0$ (cf. [Ha]).

 \skipaline \skipaline \centerline {\bf 7. A homological estimate for
 2-complexes}

 \skipaline \noindent {\bf  7.1. A norm on 1-cohomology.} 
 The author [Tu1] defined
 a homogeneous (semi-)norm on   1-cohomology of finite 2-complexes
  analogous to the
 Thurston norm on   1-cohomology of 3-manifolds.  
 We   recall this definition in a slightly generalized form. 
 By a {\sl finite 2-complex} we  mean the underlying topological space 
of a 
finite   
2-dimensional 
CW-complex   such 
that
each its point     has a neighborhood   
homeomorphic to  the cone over a finite graph. The   latter condition is aimed 
at  
eliminating  
all kinds of local wilderness. 

Let $X$ be 
   a finite
 2-complex. 
 Note that a codimension 1 subspace
  of   $X$ (in terminology of Sect. 1) is  just a  finite graph embedded 
  in $X$ and having a cylinder neighborhood in $X$. We
   call codimension 1 subspaces
  of   $X$ {\it
  regular graphs} on $X$. Any $\psi\in H^1(X)$ can be represented by
  a cooriented regular graph 
  on  $X$: it can be obtained as the preimage of a regular value of a map $X\to
  S^1$ inducing $\psi$. Set 
 $\vert\vert \psi \vert\vert_X=\min_{Y} \,  \chi_- (Y)  
$ 
where 
$Y$ 
runs 
over cooriented   regular graph 
on $ X$ such that $\psi_Y=\psi  $. 
   The   function $\psi\mapsto \vert\vert \psi \vert\vert_X:H^1(X)\to
 \Z$ is non-negative and homogeneous in the sense that $\vert\vert n\psi \vert\vert_X
 = \vert n\vert\, \vert\vert  \psi \vert\vert_X$ for any $n\in \Z, \psi \in
 H^1(X)$.
Here the inequality $\vert\vert n\psi \vert\vert_X
\leq  \vert n\vert\, \vert\vert  \psi \vert\vert_X$ is obvious (use 
$\vert n\vert$ parallel copies of any cooriented   regular graph representing $\psi$). The opposite
inequality can be obtained by the argument of [Th, p. 103] or deduced from Lemma
1.2.  We verify  in Sect. 7.4 that
 $\vert\vert \psi \vert\vert_X$ is a semi-norm. It does not depend on the choice
 of a CW-decomposition of $X$.

  \skipaline \noindent {\bf  7.2.  Theorem.  } {\sl  Let
$  K, \alpha, \Lambda= K[t^{\pm 1};\alpha] $  be as in the introduction. 
 Let $X$ be a   connected finite   2-complex  and  $\varphi : \Z[\pi_1(X)]\to \Lambda$
  be a ring homomorphism 
  compatible with a non-zero cohomology class $\psi\in  H^1(X)$ and such that
   the right 
$\Lambda$-module 
 $H_1  (X;\Lambda ) $ is finitely generated over $K$. 
Then
 $$\vert \vert \psi \vert \vert_X  \geq \max (\rk_K H_1  (X;\Lambda)-
  \delta_\varphi\,\vert \psi \vert, 0)  \eqno
 (7.1)$$
where $\vert \psi \vert$ is the maximal positive integer dividing $\psi$ 
in  $H^1(X)$ 
and $\delta_\varphi=1$ if the multiplicative group 
$\varphi(\pi_1(X))\subset \Lambda \backslash \{0\}$ 
is   cyclic and $\delta_\varphi=0$ otherwise. If   $\psi$ is 
induced   by a fibration $X\to S^1$, then $H_1  (X;\Lambda ) $ is finitely generated over $K$
    and (7.1) becomes an equality.}

  \skipaline {\sl Proof.} The proof of (7.1) goes along the same lines as the 
  proof of
  Theorem 1. First one reduces the theorem to the case of primitive $\psi$. 
  Then using Lemma 1.2 one obtains a weighted  
   cooriented regular graph  
   $Y=\cup_i (Y_i,w_i)$ on $ X$ such that 
  $X\backslash Y$ is connected, 
  $\psi_Y=\psi$, and $\vert \vert \psi \vert \vert_X=\sum_i w_i \chi_- (Y_i)$.
  If $X\backslash Y$ is not bad with respect to $\varphi$,
  then 
  by Lemma 2.3, $Y$ has no bad components so that   
  $\chi_-(Y_i)= \rk_K H_1(Y_i;K)$ for all $i$. This and Lemma 3.2 imply that 
$\vert \vert \psi \vert \vert_X  \geq \rk_K H_1  (X;\Lambda)$. 

Suppose  
that $X\backslash Y$ is   bad with respect to $\varphi$. Then $Y$ is
connected, its weight is 1,   and $\delta_\varphi=1$.  We must prove
that $ \chi_-(Y ) \geq \rk_K H_1  (X;\Lambda)-1$. Let
$F\subset K$ be a commutative subfield of $K$ defined   in Sect. 4.4. 
As there,  $H_\ast (X;\Lambda)= H_\ast (\hat X;F)\otimes_F K$
where $\hat X $ is the infinite cyclic cover  of $X$ determined by $\psi$. 
The assumptions of the theorem imply that the vector space $H_1 (\hat X;F)$ is finite
dimensional. 
Observe that 
 $Y$ lifts to a homeomorphic graph $\hat Y
\subset \hat 
X$  splitting $\hat X$ into two connected pieces, $\hat X_-$ and 
$\hat X_+$.  Let $t$ be the generating deck transformation of the covering 
$\hat 
X\to
X$ such that $t \tilde X_+ \subset \tilde X_+$. Since 
$H_1 (\hat X;F)$ is finite
dimensional,  its basis can be represented by 1-cycles lying in 
$t^{-m} \hat X_+$ for sufficiently big   $m$. Applying $t^{m}$ we obtain
   a basis of $H_1 (\hat X;F)$     represented by 1-cycles   in 
$ \hat X_+$.  Therefore 
the inclusion homomorphism  $ H_1(\hat X_+; F) \to H_1(\hat 
X, F)$ is surjective. Similarly, the inclusion homomorphism  $ H_1(\hat X_-; F) \to H_1(\hat 
X, F)$ is also surjective.
  The
Mayer-Vietoris   sequence for $\hat X=\hat X_+\cup \hat 
X_-$ 
implies
the surjectivity of the inclusion homomorphism $ H_1(\hat Y ; 
F) 
\to
H_1(\hat X; F)$.  Hence  $$\chi_- 
(Y)\geq -\chi(Y)=\dim_{F}\,
H_1(Y; F) -1 $$
$$
\geq 
\dim_{F}\,
H_1(\hat X; F)-1=\rk_K H_1  (X;\Lambda)-1 .$$ 
The proof of the last claim of the theorem follows the proof of Theorem 2 with 
obvious
changes.

 \skipaline \noindent {\bf  7.3.  Special cases.  } 
    Fix a connected finite   2-complex $X$.

 1. Let $\Lambda=F[t^{\pm
 1}]$ be the usual commutative ring of Laurent polynomials 
 on the variable $t$ with  
 coefficients in a commutative field $F$. 
 Pick $\psi \in H^1(X)$ and let $\varphi$ be the ring homomorphism $\Z[\pi_1(X)]\to \Lambda$
 defined by $\varphi(\gamma)=t^{\psi(\gamma)}$ for $\gamma\in  \pi_1(X)$. Then
 $H_1(X;\Lambda)=H_1(\hat X;F)$ where $\hat X\to X$ is the infinite cyclic
 covering determined by $\psi$.   The inequality (7.1) says that if 
 $H_1(\hat X;F)$ is finite dimensional over
 $F$ then 
 $\vert \vert \psi\vert\vert_X \geq \dim_F \,H_1(\hat X;F) -  \vert \psi\vert$.
 Using group homomorphisms $H_1(X)\to F^*$, we can obtain \lq\lq twisted" versions of
 this inequality, cf. Sect. 6.1. 
 
 2. For any $\psi \in H^1(X)$ and any $\sigma \in \Hom (\Tors H_1(X), \C^*)$, 
 it is proven in [Tu1] that 
 $$ \vert \vert \psi \vert\vert_X \geq \left \{ \begin {array} {ll}
   \vert \vert \psi \vert\vert^\sigma,~ {\rm {if}} 
\,\,\, b_1(X)\geq 2, 
\\
  \vert \vert \psi \vert\vert^\sigma- 
   \delta^1_\sigma \vert \psi\vert 
,~ 
{\rm
{if}} \,\,\,  b_1(X)=1   \end {array} 
\right.  $$
  where $\vert \vert 
...\vert 
\vert^\sigma$ is 
the Alexander-Fox  norm  on $H^1(X ) $   
determined 
by 
$\sigma$ and   
$\delta_\sigma^1=1$ if $\sigma=1$ and
$\delta_\sigma^1=0$ otherwise. 
This inequality can be deduced from Theorem 7.2
following the lines of Sect. 6.3.

    3. We can use PTFA groups in this setting exactly as in Sect. 6.4.    
   Consider a PTFA group $\Gamma$
   and 
   a  commutative field $F$. Suppose that we have   group homomorphisms
   $\varphi_0:\pi=\pi_1(X)\to \Gamma$  and $\sigma:H_1(X)\to F^*$.   
    Pick   a primitive cohomology class $ \psi\in H^1(\Gamma)$ and 
    recall the skew field $K_\psi$ and the embedding 
      $F [\Gamma]\subset K_\psi [t^{\pm 1}]$ constructed in Sect. 6.4.
       Consider the ring homomorphism 
    $\varphi: \Z[\pi] \to K_\psi [t^{\pm 1}]$ which maps any $\gamma\in \pi$ to
     $\sigma ([\gamma]) \varphi_0(\gamma) \in 
     F [\Gamma]\subset K_\psi [t^{\pm 1}]$ where $[\gamma]\in H_1(X)$ is the
     homology class
     of $\gamma$. It is clear that 
    $\varphi$ is compatible with 
     $ \varphi_0^* (\psi)\in  H^1(X)$.  
  Thus we can apply Theorem
 7.2  to estimate $\vert \vert  \varphi_0^* (\psi)\vert\vert_X$ from below. 
   
   In particular, set  $F=\C$, $ \Gamma=\pi/
   \pi^{(i)}$ with $i\geq 1$, and let $\varphi_0$ be  the   projection  
  $\pi\to \Gamma$. Let us identify $H^1(\Gamma)=H^1(\pi)$ via $\varphi_0^* $. 
  For any $\sigma\in \Hom (H_1(X),\C^*)$ and 
 any primitive $\psi \in  H^1(X)$, we obtain  the following:
       if $H_1(M; K_\psi [t^{\pm 1}])$ has finite rank over $K_\psi$, 
       then  
    $$\vert\vert \psi\vert\vert_T\geq \rk_{K_\psi} H_1(M; K_\psi [t^{\pm 1}])
    -\delta_\varphi. \eqno (7.2)$$   One can check that  the number 
     $\rk_{K_\psi} H_1(M; K_\psi [t^{\pm 1}])$  
    depends only on the restriction of $\sigma$ to $\Tors H_1(M)$. 
    The inequality (7.2) extends by linearity to
     arbitrary   $\psi \in  H^1(X)$.  For $i=1$ we recover the estimates
 discussed above. 
 
 Note that $\delta_\varphi=1$ if and only if $b_1(X)=1$,
 $\sigma (\Tors H_1(X))=1$ 
 and   either  $i=1$   or $i\geq 2$ and 
 $\pi^{(1)}=\pi^{(2)}=\pi^{(3)}=...$. In all other cases $\delta_\varphi=0$. 
  If $i\geq 2$ and 
 $\pi^{(1)}=\pi^{(2)}=\pi^{(3)}=...$ then  $H_1(M; K_\psi [t^{\pm 1}])=0$ and 
  (7.2) is of no
 interest.
    
  \skipaline \noindent {\bf  7.4.  Theorem.  } {\sl  
 For any    finite  2-complex $X$, the 
 function $\psi\mapsto \vert\vert \psi \vert\vert_X:H^1(X)\to
 \Z$ is a semi-norm.}

  \skipaline {\sl Proof.} We need to show that the function 
  $\psi\mapsto \vert\vert \psi \vert\vert_X$ is convex, i.e., $\vert\vert \psi+\psi'
 \vert\vert_X\leq 
  \vert\vert \psi \vert\vert_X+ \vert\vert \psi' \vert\vert_X$ for any
  $\psi,\psi'\in H^1(X)$.    In the case where each 1-cell  of $X$ is
  adjacent to at least two 2-cells of $X$ 
  (counting with multiplicity)  
   this is verified in [Tu1, Lemma 1.4]. The key point is that
  in this case   every vertex of
 a regular graph $Y\subset X$ is incident to at least 2 edges of $Y$
 (counting with multiplicity)
 and therefore $\chi(Y)\leq 0$ and $\chi_- (Y) =-
 \chi(Y)$. A regular graph representing  
 $\psi+\psi' $ can be
 obtained from  
 the union of regular graphs representing $\psi,\psi'$ 
 by smoothing at all crossing points. This preserves
 the Euler characteristic and hence $\vert\vert \psi+\psi'
 \vert\vert_X\leq 
  \vert\vert \psi \vert\vert_X+ \vert\vert \psi' \vert\vert_X$.
  
Consider now the general case. 
 We choose a CW-decomposition of $X$ such that
   all 2-cells  are glued to the 1-skeleton along loops  
  formed by sequences of oriented 1-cells (possibly with repetitions).
   Suppose that
  $X$ has a 1-cell $e$ which has no adjacent 2-cells. 
  All  cells of $X$ besides $e$ form
  a subcomplex $X'\subset X$. 
  If the endpoints of $e$ belong to different
  components of $X'$, then $H^1(X)=H^1(X')$. It follows from definitions and the
  fact that $\chi_- (point)=0$ that   the functions 
  $\vert\vert ... \vert\vert_X, \vert\vert ... \vert\vert_{X'}$
  coincide. If the endpoints of $e$ lie on one   
  component of $X'$, then similarly   the function  
  $\vert\vert ... \vert\vert_X$ is the composition of the restriction
  homomorphism $H^1(X)\to H^1(X')$ with $ \vert\vert ... \vert\vert_{X'}$.
  In both cases the convexity of $ \vert\vert ... \vert\vert_{X}$ would follow
  from the convexity of
  $ \vert\vert
  ... \vert\vert_{X'}$.
  
  Suppose that
  $X$ has a 1-cell $e$   adjacent  to exactly one 2-cell,
  $f$, with multiplicity 1. 
  All cells of $X$ besides $e,f$ form
  a subcomplex $X'\subset X$. It is clear that $X'$ is a deformation  retract of
  $X$ so that $H^1(X)=H^1(X')$. We claim that 
  $\vert\vert \psi \vert\vert_X= \vert\vert \psi \vert\vert_{X'}$
  for any $\psi\in H^1(X)=H^1(X')$. To see this, 
  let $e,e_1,...,e_n$ be the 1-cells of $X$ (possibly with repetitions) 
  forming the
   image of the gluing map of $f$.
 Let us  call a regular graph on $ X$ {\it normal} if
  it meets  $ f$  along  a finite set
   of  disjoint intervals   connecting   points on $e$ to points on  $\cup_i
   e_i$. Observe that for any cooriented regular graph $Z\subset X$ 
   the intersection  $ Z\cap f$ consists of a finite number of disjoint
   intervals connecting points on $e\cup \cup_i e_i$. If some of them 
    have both endpoints on $\cup_i
   e_i$, then we can push the interiors of these intervals towards $e$ and
    eventually replace each of
   them  with two intervals connecting points on  $\cup_i
   e_i$ to points on $e$. The intervals in $ Z\cap f$ with both endpoints on
   $e$ are simply eliminated. 
    These transformations produce a normal cooriented regular graph 
   $Y\subset X$
   representing the same cohomology class as $Z$. It is easy to compute from 
    the
   definition of $\chi_-$,  that $\chi_-(Y)\leq \chi_- (Z)$.
   Therefore for any $\psi\in H^1(X)$, we have $\vert\vert \psi \vert\vert_X=
   \min_Y \chi_-(Y)$ where $Y$ runs over normal cooriented regular graphs on $X$
   representing $\psi$. 
  Each such $Y$ gives rise to a cooriented regular graph $Y'=Y\cap X'$  
   on $X'$. Clearly  $\psi_{Y'}=\psi_{Y}$ and $\chi_-(Y')=\chi_-(Y)$.
   Conversely, any cooriented regular graph
   $Y'$ on $X'$    arises in this way from a
   (unique up to ambient isotopy) 
 normal  cooriented regular graph on $ X$. This implies that
   $\vert\vert \psi \vert\vert_X=  \vert\vert \psi \vert\vert_{X'}$ for any
   $\psi\in H^1(X)=H^1(X')$.  Thus, 
   the convexity of $ \vert\vert ... \vert\vert_{X}$ would follow
  from the convexity of
  $ \vert\vert
  ... \vert\vert_{X'}$.
   
   Iterating the transformations $X\mapsto X'$ as above we eventually 
    arrive to a
   2-complex $X_0$ whose all 1-cells are
  adjacent to at least two 2-cells. As we know, the function 
  $ \vert\vert ... \vert\vert_{X_0}$ is convex. Therefore the function 
  $ \vert\vert ... \vert\vert_{X}$   is also convex.

  \skipaline \noindent {\bf  7.5.  Remark.  } It follows from general
  properties of homogeneous semi-norms on lattices (cf. [Th]), that 
  the 
 function $\psi\mapsto \vert\vert \psi \vert\vert_X:H^1(X)\to
 \Z$ uniquely extends to a continuous homogeneous  semi-norm $H^1(X;\Bbb R)\to
 \Bbb R$.

   \skipaline \skipaline \centerline {\bf Bibliography}

  [Ak] M. Akaho, {\it An Estimate of Genus  of Links}, J. of Knot 
Theory and Its Ramifications, 8 (1999),
405--414.

[Au] D.  Auckly, {\it The Thurston norm and three-dimensional
Seiberg-Witten
theory}, Osaka J.  Math.  33 (1996), no.  3, 737--750.
 
 [Co]  T. Cochran, {\it Noncommutative Knot Theory},  math.GT/0206258.
 
 [COT] T. Cochran, K. Orr, P. Teichner, {\it   Knot concordance, Whitney towers
 and $L^2$ signatures}, Ann. Math., to appear.

 [Coh1] P.M. Cohn,   Algebra. Vol. 2. Second edition. 
 John Wiley and  Sons, Ltd., Chichester, 1989.  
 
  [Coh2] P.M. Cohn,  Skew fields. Theory of general division rings. Encyclopedia of Mathematics 
 and its Applications, 57. Cambridge Univ. Press, Cambridge, 1995.
 
 [Fox] R. H.  Fox, {\it Free differential calculus.  II.  The
isomorphism
problem
of
groups}, Ann.  of Math.  (2) 59 (1954), 196--210.

 [Ha] S. Harvey, {\it Higher-Order Polynomial Invariants of 3-Manifolds Giving
 Lower Bounds for the Thurston Norm}, math.GT/0207014.

[Kr] P.  Kronheimer, {\it Minimal genus in $S\sp 1\times M\sp 3$},
Invent.
Math.
135 (1999), no.  1, 45--61.

[KM1] P.  Kronheimer, T.  Mrowka, {\it The genus of embedded surfaces
in the
projective plane}, Math.  Res.  Lett.  1 (1994), no.  6, 797--808.

[KM2] P.  Kronheimer, T.  Mrowka, {\it Scalar curvature and the
Thurston norm},
Math.  Res.  Lett.  4 (1997), no.  6, 931--937.

[McM] C.  McMullen, {\it The Alexander polynomial of a 3-manifold and
the
Thur\-ston norm on cohomology}, Annales Sci. Ec. Norm. Sup. 35 (2002), 153--171.

[MST] J. Morgan, Z. Szabo, C. Taubes,   {\it 
 A product formula for the Seiberg-Witten invariants and the generalized Thom conjecture}. J. Differential
Geom. 44 (1996),  706--788.

[OS] P.  Ozsvath,  Z.  Szabo, {\it Holomorphic disks and
three-manifold
invariants:  properties and applications},  math.SG/0105202.
 
 [Pa] D. Passman,   The algebraic structure of group rings. 
 Pure and Applied Mathematics. Wiley-Interscience [John Wiley and Sons], 
 New York-London-Sydney, 1977.

[St]  B. Stenstr\"om,   Rings of quotients.
 Die Grundlehren der Math. Wissenschaften, Band 217. 
 An introduction to methods of ring theory. Springer-Verlag, New
York-Heidelberg, 1975.

[Th] W.  Thurston, {\it A norm for the homology of $3$-manifolds}, Mem.
Amer.
Math.  Soc.  59 (1986), no.  339, i--vi and 99--130.

 [Tu1] V.  Turaev, {\it A norm for the cohomology of 2-complexes},   
 Algebr. and Geom. Topology 2 (2002),   137-155 (see also math.AT/0203042). 

  [Tu2] V.  Turaev,     Torsions of 3-dimensional Manifolds.  Birkh\"auser,
Basel, to appear.
 
[Vi] S. Vidussi, {\it Norms on the cohomology of a 3-manifold and SW theory},
math.GT/0204211.

\skipaline IRMA, Univerist\'e Louis Pasteur - CNRS, 7 rue Rene Descartes, 67084
Strasbourg, France

Email:  turaev@math.u-strasbg.fr
\end{document}